\documentclass[letterpaper, 10 pt, conference]{ieeeconf}  

\usepackage[hyphens]{url}
\usepackage{hyperref}
\usepackage{breakurl}

\IEEEoverridecommandlockouts                              
\usepackage{amsmath,amssymb,multicol,bbm}
\usepackage{graphicx,color,subfigure}
\usepackage[dvipsnames*,svgnames]{xcolor}
\usepackage{mathtools}
\usepackage{etex}
\usepackage{dsfont}

\DeclareTextSymbol{\degre}{T1}{6}
\DeclareTextSymbol{\degre}{OT1}{23}

\newcommand{\ds}{\displaystyle}
\newcommand {\beq} {\begin{equation}}
\newcommand {\eeq} {\end{equation}}
\newcommand {\barr} {\begin{array}}
\newcommand {\earr} {\end{array}}
\newcommand {\bear} {\begin{eqnarray}}
\newcommand {\eear} {\end{eqnarray}}
\newcommand {\bears} {\begin{eqnarray*}}
\newcommand {\eears} {\end{eqnarray*}}
\newcommand {\bn} {\begin{aligned}}
\newcommand {\en} {\end{aligned}}

\newcommand{\beln}{\begin{multline}}
\newcommand{\eeln}{\end{multline}}
\newcommand{\eat}[1]{}

\newcommand {\mc} {\mathcal}
\newcommand {\mb} {\mathbb}
\newcommand {\ul} {\underline}
\newcommand {\ol} {\overline}

\newtheorem{prop}{Proposition}[section]
\newtheorem{problem}{Problem}[section]

\DeclareMathOperator*{\minimize}{minimize\,}
\DeclareMathOperator*{\subjectto}{subject\,to\,}

\DeclareMathOperator*{\argmin}{argmin\,}

\title{\LARGE \bf Minimizing the impact of EV charging\\ on the electricity distribution network}

\author{Olivier Beaude$^{1,2,3}$, Samson Lasaulce$^{2}$, Martin Hennebel$^{3}$, and Jamal Daafouz$^{4}$\\
$^{1}$Renault, 78280 Guyancourt, France, \\
$^{2}$L2S (CNRS -- CentraleSup\'elec-- Univ. Paris Sud 11) \\
$^{3}$GeePs -- CentraleSup\'elec, 91192 Gif-sur-Yvette, France, \\
$^{4}$CRAN, Universit\'e de Lorraine, CNRS\\
\{beaude, lasaulce\}@lss.supelec.fr, martin.hennebel@centralesupelec.fr, jamal.daafouz@univ-lorraine.fr}

\begin{document}

\maketitle

\begin{abstract}
The main objective of this paper is to design electric vehicle (EV) charging policies which minimize the impact of charging on the electricity distribution network (DN). More precisely, the considered cost function results from a linear combination of two parts: a cost with memory and a memoryless cost. In this paper, the first component is identified to be the transformer ageing while the second one corresponds to distribution Joule losses. First, we formulate the problem as a non-trivial discrete-time optimal control problem with finite time horizon. It is non-trivial because of the presence of saturation constraints and a non-quadratic cost. It turns out that the system state, which is the transformer hot-spot (HS) temperature here, can be expressed as a function of the sequence of control variables; the cost function is then seen to be convex in the control  for typical values for the model parameters. The problem of interest thus becomes a standard optimization problem. While the corresponding problem can be solved by using available numerical routines, three distributed charging policies are provided. The motivation is threefold: to decrease the computational complexity; to model the important scenario where the charging profile is chosen by the EV itself; to circumvent the allocation problem which arises with the proposed formulation. Remarkably, the performance loss induced by decentralization is verified to be small through simulations. Numerical results show the importance of the choice of the charging policies. For instance, the gain in terms of transformer lifetime can be very significant when implementing advanced charging policies instead of plug-and-charge policies. The impact of the accuracy of the non-EV demand forecasting is equally assessed.   
\end{abstract}

\vspace{-2mm}
\section{Introduction}
\label{sec:Intro}

The vast majority of current electric vehicles (EVs) charge their battery in a very simple way. The EV battery charging operation starts as soon as the user plugs his vehicle into the grid and at the maximal power which is admissible by the charging system. The merit of this charging policy is that it does not require any interaction between the user and the grid and it minimizes the time needed to reach a given charging level.  On the other hand, such a policy has the disadvantage of strongly impacting the grid since it ignores the demand profile associated with all the other devices connected to the grid; charging may typically start e.g., during the evening peak when people get back home. Since the pricing aspect is ignored by the plug-and-charge (PaC) policy, the cost paid by the user can also be affected especially in scenarios where the electricity price is time-varying. The goal of the work reported in this paper is to provide new charging policies whose main purpose is to minimize the impact of EV charging on a precise part of the electric grid namely, the distribution network (DN). More specifically, we want to minimize the impact of EV charging on the ageing of distribution transformers and on distribution Joule losses. The corresponding mathematical model can also be used for scenarios where pricing aspects are considered. The derived results can be re-exploited for other problems in smart grids such as the general problem of energy scheduling with delay constraints \cite{Bitar11,Subramanian12}; therein the system state is given by the available stored energy and the cost are market-based or generation ones.

So far, despite of the importance of the problem\footnote{In the European Union, about 5 millions distribution transformers are used and about $70\%$ of transformer failures are due to ageing and not fated events \cite{Zhang2007}. As for DN Joule losses, they represent the most important fraction of power losses in the electricity network (about two thirds in France).}, the impact of EV charging on distribution transformer ageing and Joule losses has only been addressed in a quite small number of papers. Among relevant related works we may cite \cite{Argade12,Gong2011}. The dominant approach adopted, which is well illustrated by \cite{Gong2011}, consists in exploiting a suitable model for the ageing or Joule losses, and assess the impact of charging for simple \textit{scenarios}; for instance, all EV start charging at a given time of the day (e.g., at $7$ pm) or at random instants. The \textit{algorithmic} aspect is however not developed. This is precisely what the present work proposes. This aspect of the charging problem is tackled in the related literature but, in most cases, for minimizing the monetary cost the user has to pay for recharging his vehicle (see e.g., \cite{Gan13}). Additionally, going for advanced algorithms leads to complexity issues which justify in part why considering distributed algorithms are fully relevant and even required. This observation explains why game-theoretic tools may offer the right framework to address the charging problem. In this respect, useful contributions include \cite{Agarwal2011,Mohsenian-Rad2010,Beaude12}.

The main contributions of our paper can be summarized as follows. In Sec. \ref{sec:centralizedControl}, we formulate the general (centralized) discrete-time optimal control problem to be solved to minimize the 
cost of interest namely, a linear combination of the transformer ageing and Joule losses over the DN. The system state is given by the transformer hot-spot (HS) temperature which is the most influential variable for the transformer ageing. In Sec. \ref{sec:distributed-control}, we provide three distributed charging policies which are suboptimal w.r.t. the corresponding centralized solutions but perform significantly better in terms of transformer lifetime than baseline schemes. The three proposed schemes are all based on the sequential best-response dynamics (BRD). In Sec. \ref{sec:simulations}, numerical results show the potential benefits of using the proposed charging schemes in a real system.

\section{Proposed Modeling}
\label{sec:sys_model}

In this paper, we consider a DN which comprises one transformer to which two groups of electrical devices are connected: a set of \textit{EVs} and a set of \textit{other electrical devices}. The latter is assumed to induce a power demand which is independent of
the charging policies and therefore referred to as the \textit{non-EV demand}.
The corresponding load is denoted by $\ell_t$ which is a 
deterministic function of the time and this function is always assumed to be known (except in the simulation part -Sec. \ref{sec:simulations}- where the influence of not forecasting it perfectly is assessed). Time is assumed to be slotted and indexed 
by $t \in \mathcal{T} = \{1,  \cdots, T \}$. At each time-slot $t$, of typical duration $30$ min, an EV may be active or not. The extent to which EV $i$, $i\in \mathcal{I} = \{1,\cdots, I\}$, is active on time-slot $t$ is measured by the \textit{load} it generates which is denoted by $v_{i,t}$ (this quantity will also be called the \textit{charging power} of EV $i$ at time $t$). The \textit{total load} on the transformer at time $t$ is then expressed as
\begin{equation}
\label{eq:totalLoad}
u_t= \ell_t+\sum_{i=1}^{I} v_{i,t} \ .
\end{equation}
As a useful auxiliary variable we will use the \textit{sum-EV load} at time $t$: $w_t= \sum_{i=1}^{I} v_{i,t}$ .
For the system of interest, the state is the \textit{HS temperature} which is denoted by $x_t$, $ 0 \leq  x_t \leq x_{\max}$ . A  suitable equivalent model for the evolution law of the HS temperature is as follows \cite{Rivera2008a}:
\begin{equation} 
 \label{e:theta_state} 
\forall t \in  \mc{T},  \ x_{t}= ax_{t-1}+b_1 u_t^2 + b_2 u_{t-1}^2+c_t \ , 
\end{equation}
where $ 0  \leq a  \leq 1$, $b_1 \geq 0 $, $b_2 \leq 0$, $c_t$ is a known deterministic function (it typically represents the ambient temperature in Celsius degrees, making $x_t\geq0$), and $(x_0,u_0)$ is assumed to be given.

A common and convenient way of measuring the impact of the load on the transformer ageing is to consider the instantaneous factor of accelerated ageing (FAA), which measures the speed of degradation relatively to the case of a given nominal HS temperature. Denoting the FAA by $A_t$, a well-admitted model (see \cite{IEEE1995}) is given by $A_t=e^{\alpha x_t + \beta}$, where $\alpha > 0$ and $\beta \leq 0$. For example, the case where the HS temperature is above its nominal value corresponds to $\alpha x_t + \beta >0 $ i.e., $x_t > - \frac{\beta}{\alpha}$. 
To conclude on the modeling aspect, Joule losses over time-slot $t$ are merely given by $J_{t} = K  \left(\ell_t+\sum_{i=1}^{I} v_{i,t} \right)^2$, where $K \geq 0$ is a parameter which both represents the secondary equivalent resistance of the transformer and the resistances of the different links between the transformer and the different EVs.

\vspace{-1mm}
\section{Centralized EV charging}
\label{sec:centralizedControl}

It is assumed that EV user $i$ wants the battery of his vehicle to have reached a certain state of charge (SoC) $S_i$ at time $t=T$. The corresponding constraint on the control or charging policies writes as:
\begin{equation}
\label{eq:sumV}
\forall i \in \mc{I}, \ \sum_{t=1}^T v_{i,t} \geq S_i \ . 
\end{equation}
Additionally, the charging power of EV $i$ at time $t$, $v_{i,t}$, is assumed to be non-negative and cannot exceed the maximal power at which an EV can recharge its battery: 
\begin{equation}
\label{eq:boundV}
0 \leq v_{i,t} \leq  V_{\max} \ .
\end{equation}
We denote by $x = (x_1,\cdots,x_T)$ the system \textit{state sequence} and $v = (v_{1,1},\cdots,
v_{1,T}, \cdots, v_{I,1},\cdots,v_{I,T})$  the \textit{control sequence}. The impact of the charging operation of the different EVs on the DN is measured as a composite cost which accounts for the degradation in terms of transformer lifetime and Joule losses over the whole time period  under consideration. The centralized optimal control problem of interest can then be 
formulated as follows.

\medskip

\begin{problem}[Optimal control problem formulation]\label{prob:central}  
 \begin{align}
 \displaystyle{\minimize_{v}} & \ \widetilde{C}(v,x)  =  \displaystyle{\sum_{t=1}^T e^{\alpha x_t} + f \left( \ell_t +  \sum_{i \in \mathcal{I}} v_{i,t} \right)} \label{eq:opt-control-pb} \\
 \subjectto & \ \eqref{eq:totalLoad},\eqref{e:theta_state},\eqref{eq:sumV},\eqref{eq:boundV} \text{ and } \forall t\in {\cal T}, \ x_{t} \leq x_{\max}\text{ ,} \nonumber 
 \end{align}
\end{problem}
\medskip
where $f$ is assumed to be non-decreasing and convex. It can account for effects such as Joule losses (namely, $f: s \mapsto  K s^2$) but it can also be exploited as a pricing function. Similarly, $\ell_t$ represents the non-EV demand in this paper but it can also represent the electricity fare. Note that, for ease of presentation, the scaling factor $e^{-\beta}$ is included in the function $f$ without loss of mathematical generality.

At first glance, solving problem (\ref{eq:opt-control-pb}) is a difficult task. Indeed, it is known that determining an optimal solution of an optimal control problem with saturation constraints is generally difficult, especially when the cost function is neither linear nor quadratic; here, the cost function $\widetilde{C}$ is not quadratic in the state and not necessarily quadratic in the control. However, it turns out that in the problem under investigation, the state $x_t$ can be expressed as a sole function of the sequence $(v_1, \cdots, v_{t})$  for every $t\in \mc{T}$, where $v_{t} \triangleq (v_{1,t},\cdots, v_{I,t})$.
This observation allows us to convert the initial optimal control problem into a standard optimization problem \cite{Boyd_book}. By defining the function $g_t$ as $x_t = g_t(v_1, \cdots, v_t)$ (the expression of $g_t$ is provided further), Problem \ref{prob:central} can be written as follows.

\medskip

\begin{problem}[Standard optimization prob. formulation]\label{prob:central2}  
 \begin{align}
 \displaystyle{\minimize_{v}} & \ C(v) \label{eq:standard-OP} \\
 \subjectto & \ \eqref{eq:totalLoad},\eqref{e:theta_state},\eqref{eq:sumV},\eqref{eq:boundV} \text{ and } g_t(v_1, \cdots, v_t) \leq x_{\max} \ , \nonumber
\end{align}
\end{problem}
\medskip
with $C(v)  =  \sum_{t=1}^T  e^{\alpha g_t(v_1, \cdots, v_t)} + f\left(    \ell_t +  \sum_{i \in \mathcal{I}} v_{i,t}  \right)$ .

\medskip

Note that formulating the problem as an optimization problem has a potential disadvantage. If $T$ is large, the dimension of the optimal vector(s) to be found might make any available numerical optimization routine impossible to be run, which would then necessitate to return to the initial optimal control problem formulation. For the application of interest, $T$ typically equals $24$ or $48$ if the time horizon corresponds to a day and time-slots duration is respectively an hour or half an hour. Considering up to $I=40$ EVs per distribution transformer is affordable computationally speaking. Solving the initial optimal control problem for an arbitrary $T$ appears to be an interesting direction to explore. From now on, we consider the standard optimization problem formulation. The next result can be shown.

\medskip
 
\begin{prop}\label{prop:optimProp} Optimization problem \ref{prob:central2}   
\begin{description}
\item[(i)] has at least one optimal solution if $\forall i \in \mc{I}, \ S_i \leq V_{\max} \times T$ and $\forall t \in \mc{T}, \ g_t(\tilde{v}_1, \cdots , \tilde{v}_t) \leq x_{\max}$ with $\tilde{v}_t=(\frac{S_1}{T},\cdots,\frac{S_I}{T})$;
\item[(ii)] has multiple solutions in general;
\item[(iii)] is convex if $a b_1 + b_2 \geq 0$.
\end{description}
\end{prop} 

\medskip
 
The proof is omitted. Note only that the condition in (i) is a sufficient condition to ensure that the constraint set is nonempty. Observe that the sufficient condition of (iii) means that the transformer thermal inertia (the influence of the past load or demand levels) should not be too high i.e., $-b_2 \leq a b_1$; the latter condition is satisfied for realistic values for $a$, $b_1$, and $b_2$ (see e.g.,  \cite{Rivera2008a}). To be more precise, the condition $-b_2 \leq a b_1$ is necessary and sufficient for $g_t$, which is given by
\begin{align}
& \textstyle{g_t} = \textstyle{a^t x_0 + b_1  \left(\ell_t  + \ds{\sum_{i\in \mathcal{I}}}  v_{i,t}  \right)^2 + b_2 a^{t-1} u_0^2} \label{eq:gt} \\
& \textstyle{ +(ab_1+b_2) \ds{\sum_{t'=1}^{t-1}}a^{t'-1}  \left(\ell_{t-t'} + \ds{\sum_{i\in \mathcal{I}}}  v_{i,t-t'} \right) ^2 +\ds{\sum_{t'=1}^{t}}a^{t-t'} c_{t'}} \ , \nonumber 
\end{align}
to be convex. Under this condition, $C$ is a convex function and $g_t$, $t\in\mathcal{T}$, as well. This means that the considered optimization problem is convex. This property will be directly exploited in Sec. \ref{sec:simulations} where standard convex optimization tools (Matlab function \textsf{fmincon}) are used. Some results can be provided concerning the structure of the optimal solution(s). The next proposition provides one of these results.

Let $\mc{T}_i^{\star}$ denote the set of time-slots over which EV $i$ is effectively active for a given optimal solution of Problem \ref{prob:central2}, say $v^{\star}$: $\mc{T}_i^{\star}=\lbrace t \in \mc{T}: \ v_{i,t}^{\star}>0\rbrace$. The following result holds.

\medskip

\begin{prop}
\label{prop:formOptControl} For any optimum point, we have that
\begin{description}
\item[(i)] $S_i \geq S_j   \  \Rightarrow  \  \mc{T}_j^{\star} \subset \mc{T}_i^{\star}$;
\item[(ii)] $\forall i, \ S_i =S  \ \Rightarrow \  \forall (i,j) \in \mc{I}^2, \ \mc{T}_i^{\star} = \mc{T}_j^{\star}$. 
\end{description}
\end{prop}

\medskip

The proof is omitted. A useful observation on Problem \ref{prob:central2} is that the cost function only depends on the sequence of sum-EV-loads $w=(w_1,\cdots, w_T)$. This means that the optimization problem can be solved in two steps: 1. Find an optimal sequence of sum-EV-loads; 2. Allocate the sum-EV-load among the EVs. The optimization problem associated with the determination of an optimal sequence of sum-EV-load $w$ (Step 1) is directly derived from Problem \ref{prob:central2} introducing functions $\ol{g}_t$ such that $\ol{g}_t(w_1, \cdots, w_t)=g_t(v_1,\cdots,v_t)$ and replacing constraints \eqref{eq:sumV} and \eqref{eq:boundV} respectively by  $\sum_{t=1}^T w_{t} \geq \sum_{i =1}^{I} S_i$ and $0 \leq w_t \leq I \times V_{\max}$ .

%
Since the function minimized, $\overline{C}$, is continuous and strictly convex and the inequality constraints define a convex and compact set, there is a unique solution to the sum-EV-load optimization problem. Once this problem is solved, the allocation problem associated with Step $2$ can be tackled. The latter problem is a \textit{transportation problem} where the "sources" are the $T$ time-slots with $w_t$ supply units, the "destinations" are the $I$ EVs with $S_i$ units received and $v_{i,t}$ represents the "flow" from time-slot (source) $t$ to EV $i$ (destination) \cite{Villani08}. In the case where $v_{i,t}$ is not upper bounded, there exists a feasible allocation if and only if $\sum_{i=1}^{I} S_i = \sum_{t=1}^T w_t$, which is verified here\footnote{In the sum-EV-load optimization Problem, it is easy to see that the first constraint will be saturated at optimum.}. Otherwise, finding a maximal flow in the associated graph\footnote{To be precise, a virtual source (resp. destination) has to be added and connected to each time-slot (resp. EV) with capacity $w_t$ (resp. $S_i$).} yields a feasible configuration if the value of the optimal flow obtained is $\sum_{i=1}^{I} S_i$. Possible flow search techniques will not be detailed here. More details can be found e.g., in \cite{Goldberg89}. Remarkably, the \textit{distributed} solutions we propose in the next section solve this problem by construction and transportation-theoretic tools are not necessary. Many motivations for considering distributed policies might be provided. We only mention two of them here. First, assume a scenario (called \textit{scenario $1$}) in which the control policies are computed by a single decision-making entity (e.g., an aggregator or a transformer computing device). Note that when the dimension of $v$, which is $I \times T$, becomes too high, the computational complexity for finding an optimal solution may largely exceed the available computational capacity. Therefore, even if there is one single decision-making entity, it may be required to optimize the variables of $v$ separately.  Second, another important scenario (called \textit{scenario $2$}) will be that each EV controls its own charging policy, meaning that there are $I$ controllers instead of a single one.

\section{Distributed EV charging}
\label{sec:distributed-control}

The key difference between the framework assumed in the preceding section and the present one is 
that the variables $(v_{1,1},\cdots, v_{I,1}, \cdots, v_{1,T},\cdots,v_{I,T})$ are not assumed to be controlled jointly anymore. Rather we assume they are controlled separately by $I$ decision-makers, whether the decision is taken by a single entity (scenario $1$) or effectively by the $I$ 
EVs (scenario $2$). Decision-maker $i \in \mc{I}$ therefore only controls the sequence $\ul{v}_i \triangleq (v_{i,1},\cdots, v_{i,T})$.

We propose three distributed charging policies. They are all based on the sequential best-response dynamics (BRD see e.g., \cite{lasaulce-book-2011}), which can be seen as a generalization of well-known iterative techniques such as the Gauss-Seidel method or Cournot tat\^{o}nnement.  Note that, here, we assume that the BRD algorithm is implemented \textit{offline} based on the knowledge of the sequence of non-EV load levels $(\ell_1, \cdots, \ell_T)$. Once the control policies are determined, they can be effectively run \textit{online}. In its most used form, the BRD operates sequentially such that decision-makers update their strategies in a round-robin manner. Within round $n+1$ (with $n\geq1$) the action chosen by decision-maker $i$ is computed as\footnote{If there are more than one best action, then one of them is chosen randomly.}:
\begin{equation}\label{eq:seq-BRD}
\ul{v}_i^{(n+1)} \in \argmin_{\ul{v}_i \in \mathcal{V}_i}  C\left(\ul{v}_1^{(n+1)}, \cdots,\ul{v}_{i-1}^{(n+1)}, \ul{v}_i, \ul{v}_{i+1}^{(n)}, \cdots, \ul{v}_I^{(n)} \right) 
\end{equation}

\textit{BRD for distributed dynamic charging (DDC) policies.} After Prop. \ref{prob:central2}, convexity of
$C$ w.r.t. $\ul{v}_i$ is guaranteed under the condition $ab_1 + b_2 \geq 0$, which is assumed to hold 
here. Thus, an element of the argmin set in (\ref{eq:seq-BRD}) can be obtained by solving the corresponding convex optimization problem e.g., by using known numerical techniques (e.g., using Matlab function \textsf{fmincon}). One of the assets of this distributed control policy is that complexity is reduced  compared to the centralized approach since it is linear in the number of rounds needed for convergence (say $N$, which typically equals $3$ or $4$) and the number of EVs $I$. Therefore for a numerical routine whose complexity is cubic in the problem dimension, the complexity for the centralized implementation is of the order of $I^3 T^3$ whereas it is of the order of $N I T^3$ with the distributed implementation. However, in terms of information, all the model parameters ($a$, $b_1$, $b_2$, $q$, $r$,  etc.) need to be known whether the centralized or distributed implementation is considered. If this turns out to be a critical aspect in terms of identification in practice, other techniques which only exploit directly measurable quantities such as the sum-load have to be used. This is one of the purposes of the scheme proposed next.   

\medskip

\textit{BRD and the iterative valley-filling algorithm (IVFA).} The valley-filling or water-filling charging algorithm is a quite well-known technique (see e.g., \cite{Shinwari2012}) to allocate a given additional energy need (which corresponds here to the one induced by the EVs) over time given a primary  demand profile (which corresponds here to the non-EV demand). The idea is to charge the EVs when the non-EV demand is sufficiently low. Note that this is optimal in particular when Joule losses are considered (memoryless case), i.e. $\alpha=0$. Here, the novelty relies on the fact that the proposed implementation is an iterative version of the valley-filling algorithm. Indeed, in \cite{Shinwari2012} for instance, valley-filling is used to design a scheduling algorithm but the iterative implementation is not explored. In \cite{Gan13}, a distributed algorithm which relies on a parallel implementation (the $I$ charging vectors are updated simultaneously over the algorithm iterations) is proposed. Convergence to the valley-filling solution is obtained by adding a \textit{penalty} (or stabilizing) term to the cost. Note that one of the drawbacks of the latter approach is that the weight assigned to the added term has to be tuned properly. Here, we propose a sequential version which does not have this drawback and can be seen as a power system counterpart of the iterative water-filling algorithm used in communications problems \cite{Yu02}. Convergence is ensured thanks to the exact potential property of the associated charging game (see \cite{Beaude12} for more details on the definition of this game), which is commented more at the end of the present section. At round $n+1$, the charging power of EV $i$ at time $t$ is updated as $v_{i,t}^{(n+1)} = \left[\lambda_i -\ell_t -  \sum_{j \in \mathcal{I}, j \neq i} v_{j,t}^{(n)}\right]_0^{V_{\max}}$,
where $\left[s\right]_0^{V_{\max}} =\min(V_{\max},\max(s,0))$ and $\lambda_i$ is a threshold to be chosen. The value of this threshold is obtained by setting $S_i  -  \sum_{t=1}^T v_{i,t}^{(n+1)}$ to zero\footnote{$v_{i,t}^{(n+1)}$ can be explicitly obtained in a few simple cases.}, because it is easy to see that the sum-load constraint will be active at optimum. Compared to the DDC scheme, an important practical advantage of IVFA is that it relies only on the measure of the total load $\ell_t$ (it is an "open-loop" scheme). However, both solutions are based on continuous charging power levels ($v_{i,t} \in \mb{R}$). This assumption may not be met in some real EV networks. Additionally, just as the problem of noise robustness for high-order modulations in digital communications, these two schemes may be sensitive to uncertainties on the knowledge of the non-EV demand i.e., the sequence $(\ell_1,\cdots, \ell_T)$. This motivates us to propose a third scheme, which is based on rectangular charging profiles.

\textit{BRD for rectangular charging profiles.}  The main assumption made here is that the possible strategies for the decision-makers are imposed to be rectangular charging profiles, which translates mathematically as follows:
\begin{equation}
\begin{array}{l}
\overline{\mathcal{V}}_i = \big\{ \ul{v}_i \in \mathbb{R}^T :  \forall t \in \{t_i^{\mathrm{start}}, \cdots, t_i^{\mathrm{stop}}\},  v_{i,t} = \overline{V};\\
\forall t \notin \{t_i^{\mathrm{start}}, \cdots, t_i^{\mathrm{stop}}\},  v_{i,t} = 0
\big\}
\end{array}
 \end{equation}
with $( t_i^{\mathrm{start}},   t_i^{\mathrm{stop}}) \in \mathcal{T}^2$, $ t_i^{\mathrm{start}} \leq  t_i^{\mathrm{stop}}$, and $\overline{V} \leq 
V_{\max}$ . In practice, $t_i^{\mathrm{stop}}$ may be chosen to be the minimum stopping time such that $ (t_i^{\mathrm{stop}} -    t_i^{\mathrm{start}}) \times \overline{V}  \geq S_i$. In this case, choosing the optimal charging profile amounts to choosing the optimal charging start time $t_i^{\mathrm{start}}$, which is determined by EV $i$ solving \eqref{eq:seq-BRD} "in response" to the total (except EV $i$) load sequence $(\ell_t +  \sum_{j \in \mathcal{I}, j \neq i} v_{j,t}^{(n)})_{t \in \mc{T}}$ (see \cite{Beaude12} for more details).

Motivations for using a control of this form are as follows \cite{Wu13}: 1. This strategy is easy to implement; 2. Rectangular charging profiles are believed to perform quite well in terms of battery ageing \cite{Gong2011}. From a control-theoretic point of view, also observe that a rectangular charging control can be optimal: when the state (the HS temperature here) is monotonically increasing in the control (the charging power $v$ here), it is optimal to start charging as late as possible i.e., to charge at maximal power at the end of the considered time window and charge at zero power before. However, both rectangular charging policies and IVFA charging policies are not well suited if the constraint $x_t \leq x_{\max}$ is likely to be active that is, when the maximal HS temperature of the distribution transformer can be reached. Only the DDC charging policy can easily integrate this constraint. 

To conclude this section, we provide a result which guarantees the convergence of the three proposed distributed charging policies. 

\medskip

\begin{prop}
\label{prop:convergence}[Convergence] The DDC algorithm, IVFA, and rectangular profiles-based BRD charging algorithm always converge.
\end{prop}

\medskip

This result can be proved by identifying each of the three proposed distributed policies as the sequential BRD of a certain strategic-form game. The key observation to be made is that since a common cost function (namely, $C$) is considered for the $I$ decision-makers and the individual control policies are vectors of $\mathbb{R}^T$ instead of general maps from the system state space to the charging power space, the corresponding problem can be formulated as an exact potential strategic-form game \cite{lasaulce-book-2011}. The important consequence of this is that the convergence of dynamics such as the sequential BRD is guaranteed due to the "finite improvement path" property. 
Note that although Prop. \ref{prop:convergence} provides a sufficient convergence condition for the proposed policies, characterizing the efficiency of the point(s) of convergence in comparison with the solution of Problem \ref{prob:central2} is not an easy task (study of the "Price of Anarchy" in game theory \cite{lasaulce-book-2011}), except in some special cases as presented in \cite{Beaude12}. This question will be addressed here by simulation in the following part.

%

\vspace{-3mm}
\section{Numerical analysis}
\label{sec:simulations}

The general simulation setup assumed by default is as follows. We assume that the time-slot duration is $30$ min and that an EV wants to charge its battery within a time window which starts at $5$ pm (day number $j$) to $8$ am (day number $j+1$), i.e., $T=30$; charging operations therefore take place during the evening and the night. Choosing here $f=0$, we focus here on the transformer (with memory) cost, which differentiates this contribution from the related literature often based on memoryless costs. The analysis of the simulation results in a bi-objective approach could constitute an interesting extension of this simulation part. We take $S_i=24$ kWh (capacity of a RENAULT Zoe or Fluence). During the day, we assume that the transformer load only consists of the non-EV demand $\ell_t$. We consider a $20$ kV/$410$ V transformer whose apparent power is $100$ kVA and nominal (active) power is $90$ kW (this approximately corresponds to a district of $30$ households). The transformer HS temperature evolution law is assumed to follow the ANSI/IEEE linearized Clause 7 top-oil-rise model \cite{Rivera2008a}; the corresponding parameters are $a = 0.83$, $b_1 = 30.91$, $b_2 = -19.09$, $c_t = 0.17 \times (8.47 + x^{amb}_t)$, where $x^{amb}_t$ denotes the ambient temperature at time $t$ and $x_{0}=98$\degre C (the transformer "nominal" temperature). The shut-down HS temperature is taken to be $x_{\max} = 150$\degre C. Realistic data corresponding to non-EV demand profiles and the ambient temperature are taken from the ERDF French DN Operator data basis: {\footnotesize\url{http://www.erdf.fr/ERDF_Fournisseurs_Electricite_Responsables_Equilibre_Profils}} and {\footnotesize\url{http://www.rt-batiment.fr/batiments-neufs/reglementation-thermique-2012/donnees-meteorologiques.html}}. Unless stated otherwise, the simulations are done over the $365$ days of $2012$. The transformer lifetime is inversely proportional to the average ageing: $\text{lifetime} = 40 \times T \times \left[\sum_{t=1}^T A_t\right]^{-1}$.


\begin{figure}
\includegraphics[scale=0.55]{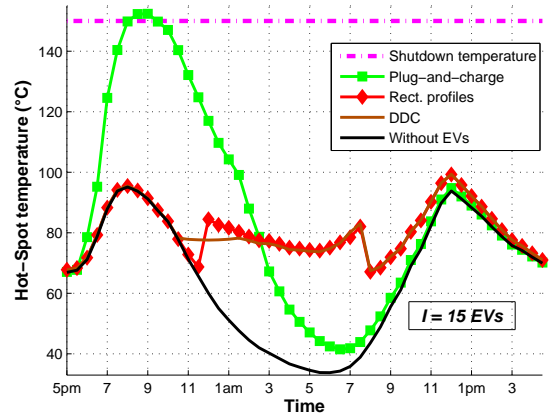}
\caption{Impact of the charging scheme on the evolution of the transformer hot-spot (HS) temperature (namely, the system state), the key variable to be controlled to manage the transformer lifetime, during the night between the $23$rd and the $24$th of March 2012. The number of EVs is here $I=15$ (penetration rate of $50$ \%). \textit{Since the instantaneous ageing is related to the HS temperature exponentially, the time-slots with the lowest temperature levels are typically preferred and heating is delayed to a large extent. Observe also that the shutdown temperature is exceeded in the plug-and-charge case, while this is never the case during the year simulated with the other charging policies (even if this was not a priori expected for the IVFA, and the strategy with rectangular profiles).}}
\label{fig:HST}
\end{figure}

The PaC policy is obtained by assuming the arrival time to follow a Poisson distribution whose mean is $\lambda=5$ ($2.5$ hours). As a reference, the scheduling policy of \cite{Shinwari2012} denoted by "SYH" according to its authors' names, will also be considered: in our case, the "hard" loads consist of the non-EV demands of each household and the "soft" loads are the EV ones; we only add the upper bound $V_{\max}$ on the "soft" load scheduled on each time-slot to be coherent with the model presented here. To assess the impact of not being able to forecast the non-EV demand perfectly, we have assumed for some figures that the optimization problems were fed with $\widetilde{\ell}_t = \ell_t + z$ where $z$ is a zero-mean additive white Gaussian noise with variance $\sigma_{\mathrm{day}}^2$. We have defined the forecasting signal-to-noise ratio (FSNR) as $\mathrm{FSNR} = 10\log_{10}\left(\frac{1}{\sigma_{\mathrm{day}}^2} \times \frac{1}{T_{\mathrm{day}}} \sum_{t=1}^{T_{\mathrm{day}}} \ell_t^2\right)$, where $T_{\mathrm{day}} = 24 \times 2$. To make the reading easier and pleasant, the figure captions have been chosen to be self-contained.
 

\begin{figure}
\includegraphics[scale=0.58]{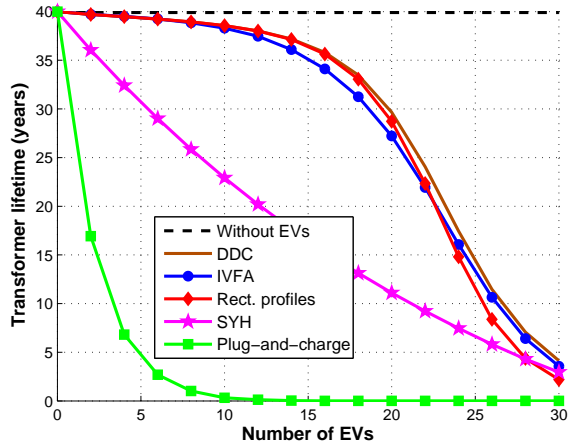}
\caption{Transformer lifetime versus the number of EVs under the assumption of perfect forecasting for the non-EV demand profile. The transformer is assumed to be chosen to be able to operate for $40$ years without EVs. \textit{The plug-and-charge policy (PaC) with Poisson arrivals is seen to be non-admissible while the proposed distributed schemes perform quite close and much better than PaC and than the policy proposed in \cite{Shinwari2012}, SYH. This latter policy schedules indeed a part of the EV loads uniformly over the time-slots, which is not suited for a cost such that the transformer one. The maximum difference between DDC and the centralized solution given by \eqref{prob:central2} is of $0.5\%$ for $I \in \lbrace 0,\cdots,30\rbrace$; this latter scenario is thus not plotted here.}}
\label{fig:lifeOfI}
\end{figure}

\begin{figure}
\includegraphics[scale=0.55]{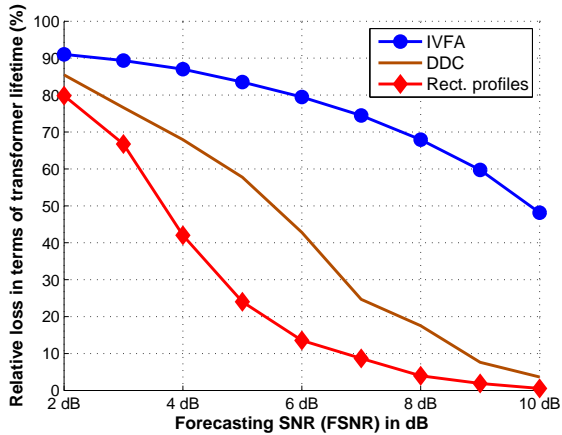}
\caption{Relative performance loss in terms of transformer lifetime versus the forecasting SNR (FSNR); the latter allows one to measure to what extent the non-EV demand can be forecasted. \textit{The most robust charging policy is the one based on rectangular profiles. Indeed, for rectangular profiles only the starting charging times need to be determined, which make them less sensitive to amplitude errors. On the other hand, the IVFA and DDC are much more sensitive to amplitude errors. Interestingly and fortunately, a typical FSNR value is $10$ dB at the scale of a district (see \cite{Shinwari2012} Fig. 10). However, with the increasing penetration of intermittent energy sources, it may be necessary to consider smaller values for the FSNR.}}
\label{fig:robust}
\end{figure}



\vspace{-1mm}
\section{Conclusion}
\label{sec:conclusion}

We have identified an important application of optimal control in the area of smart grid. The initial optimal control is difficult and left open in the case of large time horizon. As seen, for reasonable values for $I\times T$, transforming the optimal control problem into a convex optimization one is relevant since it can be solved numerically. For the cost considered, it is seen that  it is fully relevant to design distributed policies since the loss due to distributedness is typically negligible. We have also looked at the influence of the forecasting errors on the performance of these policies a posteriori. On the one hand, this points out the robustness of the simple distributed scheme with the rectangular profiles, which makes it very interesting for practical applications. On the other hand, this confirms the need to strongly integrate the forecasting aspect in the initial formulation of the problem. An interesting research direction would be then to design robust distributed dynamic charging policies. For this purpose, a stochastic formulation seems appropriate given that a good statistical knowledge can be acquired from existing power systems data bases. 

\vspace{-1mm}
\section*{Acknowledgment}
The authors would like to thank project LIMICOS - ANR-12-BS03-0005 for partly financing the work reported here.

\bibliographystyle{IEEEtran}
\bibliography{eccbib_V2}

\begin{thebibliography}{10}
\providecommand{\url}[1]{#1}
\csname url@samestyle\endcsname
\providecommand{\newblock}{\relax}
\providecommand{\bibinfo}[2]{#2}
\providecommand{\BIBentrySTDinterwordspacing}{\spaceskip=0pt\relax}
\providecommand{\BIBentryALTinterwordstretchfactor}{4}
\providecommand{\BIBentryALTinterwordspacing}{\spaceskip=\fontdimen2\font plus
\BIBentryALTinterwordstretchfactor\fontdimen3\font minus
  \fontdimen4\font\relax}
\providecommand{\BIBforeignlanguage}[2]{{%
\expandafter\ifx\csname l@#1\endcsname\relax
\typeout{** WARNING: IEEEtran.bst: No hyphenation pattern has been}%
\typeout{** loaded for the language `#1'. Using the pattern for}%
\typeout{** the default language instead.}%
\else
\language=\csname l@#1\endcsname
\fi
#2}}
\providecommand{\BIBdecl}{\relax}
\BIBdecl

\bibitem{Bitar11}
E.~Bitar, R.~Rajagopal, P.~Khargonekar, and K.~Poolla, ``The role of co-located
  storage for wind power producers in conventional electricity markets,'' in
  \emph{American Control Conference (ACC)}, 2011, pp. 3886--3891.

\bibitem{Subramanian12}
A.~Subramanian, M.~Garcia, A.~Dominguez-Garcia, D.~Callaway, K.~Poolla, and
  P.~Varaiya, ``Real-time scheduling of deferrable electric loads,'' in
  \emph{American Control Conference (ACC)}, 2012, pp. 3643--3650.

\bibitem{Zhang2007}
X.~Zhang, E.~Gockenbach, V.~Wasserberg, and H.~Borsi, ``Estimation of the
  lifetime of the electrical components in distribution networks,'' \emph{Power
  Delivery, IEEE Trans. on}, vol.~22, no.~1, pp. 515--522, 2007.

\bibitem{Argade12}
S.~Argade, V.~Aravinthan, and W.~Jewell, ``Probabilistic modeling of ev
  charging and its impact on distribution transformer loss of life,'' in
  \emph{Electric Vehicle Conference (IEVC), IEEE International}, 2012.

\bibitem{Gong2011}
Q.~Gong, S.~Midlam-Mohler, V.~Marano, and G.~Rizzoni, ``Study of {PEV} charging
  on residential distribution transformer life,'' \emph{Smart Grid, IEEE Trans.
  on}, vol.~3, no.~1, pp. 404--412, 2011.

\bibitem{Gan13}
L.~Gan, U.~Topcu, and S.~H. Low, ``{Optimal decentralized protocol for electric
  vehicle charging},'' \emph{Power Systems, IEEE Trans. on}, vol.~28, no.~2,
  pp. 940--951, 2013.

\bibitem{Agarwal2011}
T.~Agarwal and S.~Cui, ``Noncooperative games for autonomous consumer load
  balancing over smart grid,'' \emph{CoRR}, vol. abs/1104.3802, 2011.

\bibitem{Mohsenian-Rad2010}
A.-H. Mohsenian-Rad, V.~W. Wong, J.~Jatskevich, R.~Schober, and A.~Leon-Garcia,
  ``Autonomous demand side management based on game-theoretic energy
  consumption scheduling for the future smart grid,'' \emph{Smart Grid, IEEE
  Trans. on}, vol.~1, no.~3, pp. 320--331, 2010.

\bibitem{Beaude12}
O.~Beaude, S.~Lasaulce, and M.~Hennebel, ``{Charging games in networks of
  electrical vehicles},'' in \emph{Network Games, Control and Optimization
  (NetGCooP), 6th Internat. Conf. on}, 2012, pp. 96--103.

\bibitem{Rivera2008a}
L.~Rivera and D.~Tylavsky, ``Acceptability of four transformer top-oil thermal
  models: Pt. 1: Defining metrics,'' \emph{Power Delivery, IEEE Trans. on},
  vol.~23, no.~2, pp. 860--865, 2008.

\bibitem{IEEE1995}
\emph{Guide for loading mineral-oil-immersed transformers}.\hskip 1em plus
  0.5em minus 0.4em\relax IEEE Std. C57.91-1995, 1995.

\bibitem{Boyd_book}
S.~Boyd and L.~Vandenberghe, \emph{Convex Optimization}.\hskip 1em plus 0.5em
  minus 0.4em\relax Cambridge University Press, 2004.

\bibitem{Villani08}
C.~Villani, \emph{Optimal transport: old and new}.\hskip 1em plus 0.5em minus
  0.4em\relax Springer, 2008, vol. 338.

\bibitem{Goldberg89}
A.~V. Goldberg, {\'E}.~Tardos, and R.~E. Tarjan, ``Network flow algorithms,''
  DTIC Document, Tech. Rep., 1989.

\bibitem{lasaulce-book-2011}
S.~Lasaulce and H.~Tembine, \emph{{Game theory and learning for wireless
  networks: fundamentals and applications}}, Elsevier, Ed.\hskip 1em plus 0.5em
  minus 0.4em\relax Academic Press, 2011.

\bibitem{Shinwari2012}
M.~Shinwari, A.~Youssef, and W.~Hamouda, ``A water-filling based scheduling
  algorithm for the smart grid,'' \emph{Smart Grid, IEEE Trans. on}, vol.~3,
  no.~2, pp. 710--719, 2012.

\bibitem{Yu02}
W.~Yu, G.~Ginis, and J.~M. Cioffi, ``{Distributed multiuser power control for
  digital subscriber lines},'' \emph{Selected Areas in Communications, IEEE
  Journal on}, vol.~20, no.~5, pp. 1105--1115, 2002.

\bibitem{Wu13}
T.~Wu, Q.~Yang, Z.~Bao, and W.~Yan, ``{Coordinated Energy Dispatching in
  Microgrid With Wind Power Generation and Plug-in EV},'' \emph{Smart Grid,
  IEEE Trans. on}, vol.~4, no.~3, pp. 1453--1463, 2013.

\end{thebibliography}

\end{document}